\documentclass[11pt]{article}

\usepackage{amsmath, amsthm, amsfonts}
\usepackage[bottom]{footmisc}
\usepackage{geometry}
\newtheorem{lemma}{Lemma}[section]
\newtheorem{thm}{Theorem}[section]
\newtheorem{prop}{Proposition}[section]

\newtheorem{conj}{Conjecture}[section]
\newtheorem{exe}{Example}[section]
\usepackage[LY1]{fontenc}
\usepackage{graphicx}
\usepackage{float}
\usepackage{dsfont}
\usepackage{amssymb}
\usepackage{amsmath}
\usepackage{times}
\usepackage{csquotes}
\usepackage{yfonts}
\usepackage{bbm}
\usepackage{dsfont}
\usepackage{enumitem}

\title{A combinatorial proof of the Burdzy--Pitman conjecture}
\author{Stanis\l{}aw Cichomski\footnote{Faculty of Mathematics, Informatics and Mechanics, University of Warsaw (Poland). 
E-mail: s.cichomski@uw.edu.pl} , \ Fedor Petrov\footnote{St. Petersburg State University and 
St. Petersburg Department of the Steklov Mathematical Institute RAS (Russia). \newline \indent \indent E-mail: fedyapetrov@gmail.com}}

\begin{document}
\maketitle

\begin{abstract}\noindent We prove a sharp upper bound for the number of high degree differences in bipartite graphs: let  $ (U, V, E)$ be a bipartite graph with $U=\{u_1, u_2, \dots, u_n\}$ and $V=\{v_1, v_2, \dots, v_n\}$; for $n\ge k>\frac{n}{2}$ we show that $$\sum_{1\le i,j \le n}\mathbbm{1} {\Big\{|\mathrm{deg}(u_i)-\mathrm{deg}(v_j)|\ge k}\Big\}  \le  2k(n-k).$$ 
As a direct application we  show a slightly stronger,  probabilistic version of this theorem and thus confirm the Burdzy--Pitman conjecture  about the  maximal spread of coherent and independent distributions. \end{abstract}

 \section{Introduction}
 
 Let $(\Omega, \mathcal{F}, \mathbb{P})$ be a probability space. We say that a random vector $(X,Y)$ defined on this probability space, is coherent, if there exist   sub $\sigma$-fields $\mathcal{G}, \mathcal{H} \subset \mathcal{F}$ and an event $A \in \mathcal{F}$, such that
$$X   =   \mathbb{E}(\mathbbm{1}_A|\mathcal{G}), \ \ \ \ \ Y   =   \mathbb{E}(\mathbbm{1}_A|\mathcal{H}).$$
We will also say that the joint distribution of such $(X,Y)$ is coherent on  $[0,1]^2$. Hereinafter, we write $(X, Y)   \in  \mathcal{C}$ or $\mu \in \mathcal{C}$ to indicate that the vector  $(X,Y)$ or a distribution $\mu$  is coherent. Although this notation might be seen as a bit ambiguous, it does not lead to any misunderstandings.

Krzysztof Burdzy and Soumik Pal \cite{contra} prove that  for any $\delta\in (\frac{1}{2},1]$ and $(X,Y)\in \mathcal{C}$ the probability $\mathbb{P}(|X-Y|\ge \delta)$ of the difference between coherent variables exceeding a given threshold $\delta$ is bounded above by the quantity $\frac{2(1-\delta)}{2-\delta}$. 
They go on to show that this bound is sharp and it is attained by a random vector $(X,Y)$  with $X$ and $Y$ being dependent random variables.

 Let us denote
$$\mathcal{C_{I}}  =   \lbrace (X,Y): X,Y \in \mathcal{C} ,  \ X\perp Y\rbrace,$$
as a family  of those coherent distributions, which are additionally independent. In this paper we positively answer a  related and very natural question raised by   Krzysztof Burdzy  and  Jim Pitman   in \cite{pitman}, where they have formulated the following conjecture:

\begin{conj}\label{conj}
For  $\delta \in (\frac{1}{2},1]$,  we have
\begin{equation}\sup_{(X,Y)\in \mathcal{C}_{\mathcal{I}}}\mathbb{P}(|X-Y|\ge \delta)   =   2\delta(1-\delta).\label{BPC} \end{equation}
\end{conj}

\noindent Let us highlight, that this formalism should be regarded as  taking supremum over all probability spaces  $(\Omega, \mathcal{F}, \mathbb{P})$,  all events  $A\in \mathcal{F}$  and all pairs of independent sub  $\sigma$-fields  $\mathcal{G},  \mathcal{H}  \subset  \mathcal{F}$. Although there are known alternative characterizations  of coherent distributions \cite{C1, C3, C2}, let us quote \cite{pitman}:
\begin{displayquote}
For reasons we do not understand well, these general characterizations  seem to be of little help in establishing the evaluations of $\epsilon(\delta)$ [i.e. $\mathbb{P}(|X-Y|\ge \delta)$] discussed
above, or in settling a number of related problems about coherent distributions [...].
\end{displayquote}
It is our belief that this is indeed so, because of the underlying combinatorial  nature of those problems. Notice that discretization and combinatorial techniques appear already in \cite{contra, EJP}.

Let us briefly describe our approach and the organization of the paper. It is a well-known fact that the properties of two-dimensional coherent vectors are very similar to the properties of degree sequences of bipartite graphs. Remarkable example of this phenomenon can be found in \cite{tao}. Therefore, in order to  take advantage of the combinatorial nature of the problem, we start by discussing its graph-theoretic version. More precisely, in the next section we prove the following theorem.

\begin{thm} \label{graph} Let $G = (U, V, E)$ be a bipartite graph with an equal bipartition, i.e.
$$U=\{u_1, u_2, \dots, u_n\}, \ \ \ V=\{v_1, v_2, \dots, v_n\},$$
for some $n\in \mathbb{Z}_+$. For $n\ge k>\frac{n}{2}$ we have
\begin{equation} \sum_{1\le i,j \le n}\mathbbm{1} {\Big\{|\mathrm{deg}(u_i)-\mathrm{deg}(v_j)|\ge k}\Big\} \ \le \ 2k(n-k).\label{bound1} \end{equation}
\end{thm}

The proof of the  Theorem \ref{graph} is based on an idea similar to the spread bounding theorem of Erd\H{o}s et al. -- see  \cite{erdos}. In Section $2$ we then provide an elementary example showing that the  bound (\ref{bound1}) is sharp.
In what follows in Section $3$, we show how to reduce the initial problem to the Theorem \ref{graph}. To this end, we  make use of an appropriate  sampling construction,  similar in spirit to \cite{lovasz}. The key idea is to approximate  a fixed coherent distribution with a randomly  generated sequence of graphs. We then apply Theorem \ref{graph} to each of the graphs in the sequence and obtain (\ref{BPC})  by passing to the limit.
 
\section{Number of high degree differences in bipartite graphs}

Let $G = (U, V, E)$ be a bipartite graph with an equal bipartition, that is a triplet
$$U=\{u_1, u_2, \dots, u_n\}, \ \ \ V=\{v_1, v_2, \dots, v_n\},$$
and
$$E\subset U\times V,$$
for some fixed $n\in \mathbb{Z}_+$. Let us also choose a natural number $k$ satisfying $n\ge k>\frac{n}{2}$. Hereinafter, we denote the degree sequences of $G$ as $(\alpha_i)_{i=1}^{n}$ and $(\beta_j)_{j=1}^n$, i.e. $\alpha_i =\mathrm{deg}(u_i)$ and $\beta_j=\mathrm{deg}(v_j)$ for all $1\le i, j \le n$. Without loss of generality we also assume that 
$$\alpha_1\ge \alpha_2\ge\dots\ge \alpha_n,$$
$$\beta_1 \ge \beta_2\ge \dots \ge \beta_n.$$

We start with an observation similar to the spread bounding theorem of Erd\H{o}s et al. -- see  \cite{erdos}.

\begin{lemma} \label{st} There exist $s, t \in \{1,2,\dots, n-k+1\}$ such that $\alpha_s \le \beta_{s+k-1}+k-1$ and $\beta_t \le \alpha_{t+k-1}+k-1$. \end{lemma}
\noindent \textit{Proof}: We will prove only the existence of $s$, as the case of $t$ is analogous. Assume for the sake of contradiction that such a number $s$ does not exists. Therefore, the total number of edges incident to $u_1,u_2,\dots, u_{n-k+1}$ is at least $\beta_k+\beta_{k+1}+\dots+\beta_n+k(n-k+1)$. Observe that  at least $k(n-k+1)$  of these edges go to  vertices $v_1,v_2,\dots, v_{k-1}$. Let us denote
$$\tilde{E}:=E\cap \Big(\{u_1,u_2,\dots, u_{n-k+1}\}\times \{v_1, v_2,\dots, v_{k-1} \} \Big),$$ and notice that we have just shown that $|\tilde{E}|\ge k(n-k+1)$. On the other hand, we clearly have $$|\tilde{E}|\le (k-1)(n-k+1),$$
which is a contradiction. \qed
\\ 

We now prove the Theorem \ref{graph}. With this result, we  establish a natural upper bound on the number of possible pairs of vertices with high degree differences in $G$. \\ \\ \textbf{Proof of Theorem \ref{graph}}: For $1\le i,j \le n$, let us call $(i,j)$ an  $\mathcal{A}$-pair if $\alpha_i\ge \beta_j+k$. Correspondingly, let us call $(i,j)$ a $\mathcal{B}$-pair if $\beta_j\ge \alpha_i+k$. Since  $k>\frac{n}{2}$, we have $\alpha_i>\frac{n}{2}$ for all $\mathcal{A}$-pairs $(i,j)$ and $\alpha_i<\frac{n}{2}$ for all $\mathcal{B}$-pairs $(i,j)$. As a consequence, there exists an $i_0\in\{1,2,\dots, n+1\}$ such that:
\begin{enumerate}
\item $i\le i_0-1$ for any $\mathcal{A}$-pair $(i,j)$,
\item $i\ge i_0$ for any $\mathcal{B}$-pair $(i,j)$.
\end{enumerate} Analogously, there exists $j_0\in\{1,2,\dots, n+1\}$ such that:
\begin{enumerate}  \setcounter{enumi}{2}
\item $j\le j_0-1$ for any $\mathcal{B}$-pair $(i,j)$,
\item $j\ge j_0$ for any $\mathcal{A}$-pair $(i,j)$.
\end{enumerate}Observe that by the Lemma \ref{st}, we also have:
\begin{enumerate}  \setcounter{enumi}{4}
\item for any $\mathcal{A}$-pair $(i,j)$ either $i<s$ or $j>s+k-1$, 
\item for any $\mathcal{B}$-pair pair $(i,j)$ either $j<t$ or $i>t+k-1$.
\end{enumerate}

We now show that the restrictions $1$--$6$, regardless of the initial graph, imply that the total number of $\mathcal{A}$-pairs and $\mathcal{B}$-pairs together is at most $2k(n-k)$. Let us fix $i_0, j_0 \in \{1,2,\dots, n+1\}$. First, we verify that the optimal values of $s$ and $t$ satisfy $s,t \in \{1, n-k+1 \}$. 
\\

Notice that the variable $s$ appears only in the $5$-th condition and thus the value of $s$ is not important for bounding the number of $\mathcal{B}$-pairs. Moreover, observe that if $i_0\le n-k+1$, then for $s=n-k+1$ the condition $5$. is automatically fulfilled and thus $s=n-k+1$ is an optimal value. Similarly, if $j_0\ge k+1$, then for $s=1$ the condition $5$. is also automatically fulfilled and $s=1$ is an optimal value. Finally, let us assume that  $i_0\ge n-k+2$ and $j_0\le k$. In this case, the restrictions imposed by the condition $5$. remove exactly $(i_0-s)(s+k-j_0)$ additional pairs. Therefore, as the last expression is a concave function of $s\in[1,n-k+1]$, it is minimized in one of the endpoints.  Hence we may assume that $s=1$ or $s=n-k+1$, as desired.  Analogously, we show that $t=1$ or $t=n-k+1$ is optimal. There are four possible cases now:
\begin{enumerate}[label=\alph*.]
\item $s=1, \ t=n-k+1$. We have $j\ge k+1$ for all $\mathcal{A}$-pairs and $j\le n-k$ for all $\mathcal{B}$-pairs $(i,j)$. Thus any $i$ participates in at most $n-k$ of $\mathcal{A}$-pairs and in at most $n-k$ of $\mathcal{B}$-pairs. Therefore,  since a fixed vertex can not participate in both types of pairs,  every $i$ participates overall in at most $n-k$ pairs. As a consequence, the total number of pairs does not exceed $n(n-k)<2k(n-k)$.
\item $s=n-k+1, \ t=1$. This case is symmetric to the previous one.
\item $s=1, \ t=1$.  We have $j\ge k+1$ for all $\mathcal{A}$-pairs and $i\ge k+1$ for all $\mathcal{B}$-pairs $(i,j)$. Let us denote $a:=\max(k+1,j_0)$ and $b:=\max(k+1,i_0)$. Then the total number of $\mathcal{A}$-pairs is bounded by $(n-a+1)(b-1)$, while the total number  of $\mathcal{B}$-pairs is at most $(n-b+1)(a-1)$. Notice, that for $a,b\in [k+1,n+1]$ the sum 
$$S:=(n-a+1)(b-1)+(n-b+1)(a-1),$$ 
is bilinear and it is maximized at one of four endpoints. For $a=b=k+1$, we get $S=2k(n-k)$. For, say $a=n+1$, we get $S=n(n-b+1)\leqslant n(n-k)< 2k(n-k)$.
\item $s=n-k+1, \ t=n-k+1$. This case is analogous to c.
\end{enumerate}
Hence we have shown that  the Theorem \ref{graph} holds in all cases. This ends the proof. \qed
\\

We end this section with an example showing that the  upper bound $2k(n-k)$ in (\ref{bound1}) cannot be improved. Note that a straightforward modification of this example shows that $2\delta(1-\delta)$ in (\ref{BPC})  is also sharp.

\begin{exe}Consider $n,k\in\mathbb{Z}_+$, with $n\ge k >\frac{n}{2}$. Let $G_{n,k}=(U, V, E)$, where $U=V=\{1,2,\dots, n\}$ and
$$E \ = \ \{(u,v)\in U\times V : \ \max(u,v)\le k\}.$$
We clearly have $$ \sum_{1\le i,j \le n}\mathbbm{1} {\Big\{|\mathrm{deg}(u_i)-\mathrm{deg}(v_j)|\ge k}\Big\} \ = \ 2k(n-k).$$ \end{exe}
Moreover, one can  check that inequality (\ref{bound1}) becomes an equality exactly for those graphs $G$ that are isomorphic to $G_{n,k}$ or to its complement $\overline{G}_{n,k}$. This follows  easily from the proof of Theorem \ref{graph} and we leave the details to interested reader.

\section{Proof of the Burdzy--Pitman conjecture}

By $\mathcal{C}_\mathcal{I}(n)$ we denote the set of those $(X, Y) \in \mathcal{C}_{\mathcal{I}}$, that both $X$ and $Y$ take at most $n$ different
values. 

\begin{prop}\label{xnyn} Let $(X, Y)$ be coherent and independent, and let $n$ be a positive integer. Then there exists $(X_{n},Y_{n})\in \mathcal{C}_{\mathcal{I}}(n)$, such that  $|X-X_{n}|  \le  \frac{1}{n}$ and $|Y-Y_{n}|  \le  \frac{1}{n}$, almost surely.   \end{prop}
\noindent The proof of the above Proposition can be found in \cite{mastersthesis, contra}. In what follows, fix any $\delta \in (\frac{1}{2},1]$.
 
\begin{prop} \label{finite} To prove the Conjecture \ref{conj} it is enough to verify it for $ \ \bigcup_{n=1}^{\infty}\mathcal{C}_{\mathcal{I}}(n)$. \end{prop}
\noindent \textit{Proof}: Fix $(X,Y)\in \mathcal{C}_{\mathcal{I}}$ and choose $(X_n, Y_n)$ as in Proposition \ref{xnyn}. By the triangle inequality we get
$$\mathbb{P}(|X-Y|\ge \delta) \ \le \ \mathbb{P}(|X_n-Y_n|\ge \delta-2/n).$$
Thus, assuming that the Conjecture \ref{conj} is true for $ \ \bigcup_{n=1}^{\infty}\mathcal{C}_{\mathcal{I}}(n)$, for $n$ large enough we obtain
$$\mathbb{P}(|X-Y|\ge \delta) \ \le \ 2(\delta-2/n)(1-\delta+2/n).$$
Passing to the limit ends the proof. \qed
\\ 

We are now able to prove our main result. \\ \\ \textbf{Proof of Conjecture \ref{conj}}: Fix $(X,Y)\in \bigcup_{n=1}^{\infty}\mathcal{C}_{\mathcal{I}}(n)$. There exists a probability space $(\Omega, \mathcal{F}, \mathbb{P})$, independent sub $\sigma$-fields $\mathcal{G}, \mathcal{H} \subset \mathcal{F}$ and an event $A \in \mathcal{F}$, such that $X = \mathbb{E}(\mathbbm{1}_A|\mathcal{G})$ and $Y=\mathbb{E}(\mathbbm{1}_A|\mathcal{H})$. Furthermore, for some $N,M\in \mathbb{Z}_+$, we may suppose that $X$ takes values $x_1, x_2, \dots, x_N$ on sets $G_1,G_2,\dots, G_N$, respectively of the probability space $(\Omega, \mathcal{F}, \mathbb{P})$, and $Y$ takes values $y_1,y_2,\dots, y_M$ on sets $H_1,H_2,\dots, H_M$. For simplicity, we can also assume that
$$\mathcal{G}= \sigma \Big(G_1, G_2, \dots, G_N \Big),$$
$$\mathcal{H}=\sigma \Big(H_1, H_2, \dots, H_M \Big),$$
meaning that $\sigma$-fields $\mathcal{G}, \mathcal{H}$ are generated by those disjoint partitions of $\Omega$.
For $1 \le i \le N$ and $1 \le j \le M$,  denote the probabilities $p_i = \mathbb{P}(G_i), q_j = \mathbb{P}(H_j)$ and 
$$\rho_{i,j} = \frac{\mathbb{P}(G_i\cap H_j \cap A)}{\mathbb{P}(G_i\cap H_j)}.$$
Then by the independence we have $\mathbb{P}(G_i\cap H_j)=p_iq_j$ and
 \begin{equation} x_i =\sum_{j=1}^{M}q_j\rho_{i,j}, \ \ \ \ \ \ 1 \le i \le N, \label{sumxi} \end{equation}
 \begin{equation}  y_j=\sum_{i=1}^{N}p_i  \rho_{i,j},  \ \ \ \ \ \ 1 \le j \le M, \label{sumyj} \end{equation}
 which follows from a direct computation.
\\

First, we show how to construct a sequence of bipartite graphs $G_n=(U_n, V_n, E_n)$ with $|U_n|=|V_n|=n$, such that:
\begin{enumerate}
\item in $U_n$ we have \ $p_in+O(n^{3/4})$ \ vertices of degree \ $x_in+O(n^{3/4})$, \ $i=1,2,\dots, N$,
\item in $V_n$ we have \ $q_jn+O(n^{3/4})$ \ vertices of degree \ $y_jn+O(n^{3/4})$, \ $j=1,2,\dots, M$,
\end{enumerate}
where by $O(n^{3/4})$ we denote any quantity bounded in magnitude by $Cn^{3/4}$ for some absolute constant $C > 0$, which is
uniform in $i$ and $j$. 


To this end, let us fix $n$ and without loss of generality assume that $n$ is large. We choose $n$ independent points $u_1, u_2, \dots, u_n$ in our  probability space $(\Omega, \mathcal{F}, \mathbb{P})$ and for $1\le i \le n$ denote $\alpha_i=s$ if $u_i\in G_s$. In other words, $(\alpha_1, \alpha_2, \dots, \alpha_n)$ is an i.i.d. sample from the set $\{1,2,\dots, N\}$ with weights $p_1, p_2, \dots, p_N$, respectively. We can think about this sample as a randomly generated sequence of labels. Let $A_s=\sum_{i=1}^{n}\mathbbm{1}_{\{\alpha_i=s\}}$ be the number of labels equal to $s$, $1\le s\le N$.
Observe that $A_s$ is clearly a sum of $n$ independent Bernoulli random variables. Hence, by the well known Hoeffding's inequality, we have 
$$\mathbb{P}(|A_s- n p_s|\ge nr) \ \le \ 2\cdot e^{-2nr^2},$$
for all positive $r$. Consequently, setting $r=n^{-1/4}$ we get
$$\mathbb{P}(|A_s- n p_s|\ge n^{3/4}) \ \le \ 2\cdot e^{-2\sqrt{n}}.$$
Thus, as $n$ is large, with high probability we have $A_s=np_s+O(n^{3/4})$ for all $1\le s \le N$. In fact, if this would not be the case, we can always reject the labels $(\alpha_1, \alpha_2, \dots, \alpha_n)$ and resample them again.
 Analogously, we choose points $v_1, v_2, \dots, v_n$ and generate an i.i.d. sample $(\beta_1, \beta_2, \dots, \beta_n)$ from the set $\{1,2,\dots, M\}$ with weights $q_1,q_2,\dots, q_M$, respectively. Next, let $B_t=\sum_{j=1}^{n}\mathbbm{1}_{\{\beta_j=t\}}$ be the number of labels equal to $t$, $1\le t\le M$. As previously, we can further assume that  $B_t=nq_t+O(n^{3/4})$ for all $t$.
 \\
 
Next, conditioned on the labels $(\alpha_i)_{i=1}^{n}$ and $(\beta_j)_{j=1}^{n}$, we generate a random bipartite grah $(U_n, V_n, E_n)$:
\begin{enumerate}
\item for $1\le i,j\le n$ independently, generate an idicator variable $Z_{i,j}$ with
$$ \mathbb{P}_{\alpha, \beta}(Z_{i,j}=1)  = 1-\mathbb{P}_{\alpha, \beta}(Z_{i,j}=0)  =  \rho_{\alpha_i, \beta_j},$$
where $(\Omega_{\alpha, \beta}, \mathcal{F}_{\alpha, \beta}, \mathbb{P}_{\alpha,\beta})$ is a new  space on which such independent $(Z_{i,j})_{i,j}$ can be constructed.
\item for $1\le i,j\le n$, \ set $(u_i,v_j)\in E_n$ \ iff \ $Z_{i,j}=1$, 
\end{enumerate} 
or equivalently, let $Z_{i,j}=\mathbbm{1}_{\{(u_i,v_j) \in E_n\}}$. For $1\le i \le n$, we can now write
$$\mathbb{E}_{\alpha, \beta}\ \mathrm{deg}(u_i) \ = \ \mathbb{E}_{\alpha, \beta}\Bigg(\sum_{j=1}^{n}Z_{i,j} \Bigg) \ = \ \sum_{t=1}^{M}B_t\rho_{\alpha_i, t}\ = \ \sum_{t=1}^{M}\Big(nq_t+O(n^{3/4})\Big) \rho_{\alpha_i, t},$$
and hence
 \begin{equation} \mathbb{E}_{\alpha, \beta}\ \mathrm{deg}(u_i) \ = \ nx_{\alpha_i}+O(n^{3/4}), \label{Exi} \end{equation}
where the last line follows from (\ref{sumxi}). Similarilly, for $1\le j \le n$, by (\ref{sumyj}) we get 
\begin{equation} \mathbb{E}_{\alpha, \beta}\ \mathrm{deg}(v_j) \  = \ ny_{\beta_j}+O(n^{3/4}). \label{Eyj}\end{equation}
Lastly, again by the Hoeffdings's inequality, we have
\begin{equation}\mathbb{P}_{\alpha, \beta}\Big(|\mathrm{deg}(u_i)-\mathbb{E}_{\alpha, \beta} \ \mathrm{deg}(u_i)| \ge n^{3/4}\Big) \ \le \ 2\cdot e^{-2\sqrt{n}}, \label{H1} \end{equation}
and
\begin{equation}\mathbb{P}_{\alpha, \beta}\Big(|\mathrm{deg}(v_j)-\mathbb{E}_{\alpha, \beta} \ \mathrm{deg}(v_j)| \ge n^{3/4}\Big) \ \le \ 2\cdot e^{-2\sqrt{n}}, \label{H2} \end{equation}
for all $i,j \in \{1,2,\dots, n\}$. Note that the concentration rates (\ref{H1}) and (\ref{H2}) are exponential in $\sqrt{n}$. Thus, since $n$ is large, with high probability all these concentrations take place. Then, by (\ref{Exi}) and (\ref{Eyj}), we have $ \mathrm{deg}(u_i) =  nx_{\alpha_i}+O(n^{3/4})$ and $ \mathrm{deg}(v_j)  =  ny_{\beta_j}+O(n^{3/4})$ for all $i,j\in\{1,2,\dots, n\}$. This, together with bounds on $(A_s)_{s=1}^{N}$ and  $(B_t)_{t=1}^{M}$, proves  that $G_n$  does indeed satisfy the structural conditions stated.\\

In what follows, we add additional subscripts and write $u_{i}^{(n)}$ and $v_{j}^{(n)}$ for generic elements of $U_n$ and $V_n$, respectively. We can now write 
$$\mathbb{P}(|X-Y|\ge \delta) \ = \ \sum_{\substack{1\le i \le N \\ 1\le j \le M}}\mathbbm{1}_{\{ |x_i-y_j|\ge \delta \}} \cdot p_iq_j $$
\begin{equation} = \ \lim_{n\to \infty} \ \frac{1}{n^2} \sum_{\substack{1\le i \le N \\ 1\le j \le M}}\mathbbm{1}_{\{ |nx_i-ny_j|\ge n\delta \}} \cdot \Big(p_in+O(n^{3/4})\Big)\Big(q_jn+O(n^{3/4})\Big).\label{limit1}\end{equation}
Observe, that by the  triangle inequality and defining properties of $G_n=(U_n, V_n, E_n)$, we have
$$|nx_{\alpha_i}-ny_{\beta_j}| \ \le  \ |\mathrm{deg}(u_i^{(n)})-\mathrm{deg}(v_j^{(n)})|+2\cdot O(n^{3/4}),$$
for all $i,j\in\{1,2,\dots, n\}$. Thus, we can further estimate the bound (\ref{limit1}) by
$$\le \ \limsup_{n\to \infty} \ \frac{1}{n^2} \sum_{1\le i, j \le n}\mathbbm{1}\Big\{ |\mathrm{deg}(u_i^{(n)})-\mathrm{deg}(v_j^{(n)})|\ge n\delta-2 O(n^{3/4}) \Big\}.$$
Finally,  applying Theorem \ref{graph} to each of the bipartite graphs $G_n$, we obtain
$$\le  \ \limsup_{n\to \infty}\ \frac{1}{n^2} \cdot 2 \Big(n\delta -2 O(n^{3/4})\Big)\Big(n-n\delta+2 O(n^{3/4})\Big) \ = \ 2\delta(1-\delta),$$
which ends the proof. \qed

\paragraph*{Acknowledgments.} The work of the second-named author is supported by  \emph{Junior Leader} grant  of \newline Theoretical Physics and Mathematics Advancement Foundation ``BASIS''.


\setlength{\baselineskip}{2ex}

\bibliographystyle{plain}
\bibliography{BiPv2}

\end{document}